\documentclass[12pt]{article}
\title{URYSOHN LEMMAS IN TOPOLOGICAL VECTOR SPACES}
\usepackage{amsmath}
\author{{ Ramkumar. S and Ganesa Moorthy. C}\\
{\small Department of Mathematics, Alagappa University, Karaikudi - 630 003, India.}\\
 {\small Email: $1$ ramkumarsolai@gmail.com\ \ and    $2$ ganesamoorthyc@gmail.com}}
\date{}
\begin{document}

\maketitle
\newtheorem{lem}{\sf Lemma}
\newtheorem{defn}[lem]{\sf Definition}
\newtheorem{thm}[lem]{\sf Theorem}
\newtheorem{fir}[lem]{\sf First Urysohn lemma}
\newtheorem{pro}[lem]{\sf Proposition}
\newtheorem{snd}[lem]{\sf Second Urysohn lemma}

\begin{abstract}
Two variations of classical Urysohn lemma for subsets of topological vector spaces are obtained in this article. The continuous functions constructed in these lemmas are of quasi-convex type.
\end{abstract}
\section*{Introduction}

The classical technique for the proof of Urysohn's lemma (See: \cite{Jam}) is applied to derive metrization theorem (See: \cite{Rud}) for topological vector spaces. It is possible to find the continuous functions of quasi-convex type, when the same technique is applied on certain subsets of  a topological vector spaces. The Urysohn lemma is derived on normal spaces or on locally compact spaces. So, two concepts of convex normal and convex regular subsets of a topological vector spaces are introduced; their properties are studied; and two Urysohn lemmas are derived in this article.\\
\section*{Main Results}

\begin{defn}

A subset $A$ of a topological vector space is said to be locally convex if every point has a local base in $A$ consisting of convex subsets of $A$.
\end{defn}
\begin{defn}

Let $A$ be a locally convex subset of a topological vector space $X$. A set $A$ is said to be convex regular, if for a given point $x$ and a given open convex set $B$(open in $A$) such that $x\in B$, there is an open convex subset $C$ of $A$ such that $x\in C\subseteq cl\ C\subseteq B$.
\end{defn}
\begin{defn}

Let $A$ be a locally convex subset of a topological vector space $X$. A set $A$ is said to be convex normal, if for a given non empty closed set $A_0$ and a given open convex set $B$ such that $A_0\subseteq B$, there is an open convex subset $C$ of $A$ such that $A_0\subseteq C\subseteq cl\ C\subseteq B$.
\end{defn}
\begin{pro}

Every compact convex set in a locally convex topological vector space is convex normal.
\end{pro}
{\textsf Proof:}
Let $A_0$ be a closed convex subset of a compact convex set $A$ in a locally convex topological vector space $X$. Let $B$ be an open convex set in $A$ containing $A_0$. To each $x\in A_0$, find an open convex neighbourhood $U_x$ of 0 in $X$ such that $(x+U_x)\mathop{\cap}A\subseteq (x+cl\ (U_x))\mathop{\cap}A \subseteq B$. Since $A_0$ is compact, find a finite subset $\{x_1,x_2\cdots x_n\}$ of $A_0$ such that
\begin{equation} A_0\subseteq \mathop{\cup}\limits_{i=1}^n (x_i+U_{x_i})\mathop{\cap}A \label{e1} \end{equation} The convex hull $C$ of right hand side of (\ref {e1}) is open in $A$ and it is contained in $B$. Note that, any element in $C$ is of the form $\lambda_1y_1+\lambda_2y_2\cdots \lambda_ny_n$ with $y_i\in x_i+U_{x_i},\ 0\leq \lambda \leq 1$ and $\sum\limits_{i=1}^n \lambda_i =1$. For each $y\in cl\ C$, there is net  $\left(\lambda_{1\delta} y_{1\delta} +\lambda_{2\delta}y_{2\delta}\cdots \lambda_{n\delta}y_{n\delta}\right)_{\delta \in \mathcal{D}}$ in $C$ which converges to $y$. We can find a sub net $\left(\lambda_{1\alpha} y_{1\alpha} +\lambda_{2\alpha}y_{2\alpha}\cdots \lambda_{n\alpha}y_{n\alpha}\right)_{\alpha \in \mathcal{G}}$ of this net such that $\lambda_{1\alpha}\to\lambda_1,\ y_{1\alpha}\to y_1,\cdots \lambda_{n\alpha}\to \lambda_n, y_{n\alpha}\to y_n$(say). Since $0\leq \lambda \leq 1$, $\sum\limits_{i=1}^n \lambda_i =1$, $y_i\in x_i+cl\ (U_{x_i})$ and the convex hull of $\mathop{\cup}\limits_{i=1}^n (x_i+cl\ (U_{x_i}))\mathop{\cap}A$ is contained in $B$, we conclude that $y=\sum\limits_{i=1}^n\lambda_iy_i\in B$. Thus $cl\ C\subseteq B$. This proves the result.

\begin{defn}
A function $f$ from a convex subset $A$ of real vector space into the real line is said to be 
\begin{enumerate}
\item[(i)] convex if $f(\lambda x+(1-\lambda)y) \leq \lambda f(x)+(1-\lambda)f(y)$, for every $x,y \in A$, for every $\lambda\in [0,1]$.
\item[(ii)] quasi convex if $\{x\in A:f(x)<r\}$(or equivalently $\{x\in A:f(x)\leq r\}$) is a convex set,for every real number $r$. 
\end{enumerate}
\end{defn}
\begin{defn}
Let $A$ be a convex subset of a topological vector space $X$. The set $A$ is said to be `(quasi) convex completely regular' if for a given point $x\in A$ and a given open convex set $B\subseteq A$ such that $x\in B$, there is continuous (quasi)convex function $f:A\to [0,1]$ by $f(x)=0$ and $f(A\backslash B)=\{1\}=f(B^c)$.
\end{defn}
\begin{pro}
Every convex subset of a locally convex space is convex completely regular, when it contains more than one point.
\end{pro}
{\textsf Proof:} Let $A$ be a locally convex subset of a locally convex space. Let $x\in A$ and $B$ be an open convex subset of $A$ such that $x\in B$. Without loss of generality, let us assume that there is a continuous semi norm $p$ on $X$ such that $\inf\{p(x-y):y\in A\backslash B\}\geq 1$. Define a map $f:A\to [0,1]$ by $f(y)=\min\{p(x-y),\ 1\}$. Then $f$ is continuous convex mapping such that $f(x)=0$ and $f(y)=1$, for every $y\in A\backslash B$.
\begin{fir}
Let $A$ be a convex normal subset of a topological vector space $X$. Let $A_0$ be a non empty closed convex subset of $A$ and $B$ be an open convex subset of $A$ containing $A_0$. Then there is a continuous quasi convex function  $f:A\to [0,1]$ such that $f(A_0)=\{0\}$ and $f(A\backslash B)=\{1\}$.
\end{fir}
{\textsf Proof:} Let $P$ be the set of all rational numbers in the interval [0,1]. Define $U_1=B$. Find an open convex subset $U_0$ in $A$ such that $A\subseteq U_0\subseteq cl (U_0)\subseteq U_1 \subseteq B$. Find $(U_r)_{r\in P}$ of open convex sets such that $` p<q \Rightarrow cl(U_p)\subseteq U_q$'. Define $f:A\to [0,1]$ by $f(x)=\begin{cases} \inf\ \{p\in P:x\in U_p\},&\mbox{  if  } x\in U_1\\ 1&\mbox{  if  } x\notin U_1 \end{cases}$

By the proof of the classical Urysohn lemma (See:\cite[Theorem 33.1]{Jam}), $f$ is a continuous function. Since $\{x\in A:f(x)\leq\epsilon\}$ is convex, for every $\epsilon \geq 0$, $f$ is a quasi convex function.
\begin{lem}
Let $A$ be a locally convex locally compact subset subset of a locally convex space $(X,(p_i)_{i\in I})$. Suppose further that $A$ is complete as a uniform space with the uniformity induced by $(p_i)_{i\in I}$. Let $A_0$ be a nonempty compact convex subset of $A$ and $B$ be an open(in $A$) convex subset of $A$ such that $A_0\subseteq B$. Then there is an open convex subset $C$ of $A$ such that $A_0\subseteq C \subseteq cl C\subseteq B$ and $cl C$ is compact.
\end{lem}
{\textsf Proof:} To each $x \in A_0$, find an open convex neighbourhood $U_x$ of 0 in $X$ such that $(x+U_x)\mathop{\cap}A\subseteq (x+cl(U_x))\mathop{\cap}A \subseteq B$, and $(x+cl(U_x))\mathop{\cap}A $ is compact. Find a finite subset $\{x_1,x_2\cdots x_n\}$ of $A_0$ such that \begin{equation} A_0\subseteq \mathop{\cup}\limits_{i=1}^n (x_i+U_{x_i})\mathop{\cap}A \label{e2} \end{equation} The convex hull $C$ of right hand side of (\ref {e2}) is open in $A$ and it is contained in $B$. By theorem in \cite{Rud} and by completeness of $A$, the convex hull of 
$\mathop{\cup}\limits_{i=1}^n (x_i+cl(U_{x_i}))\mathop{\cap}A $ has a compact closure. Thus $cl C$ is contained in a compact subset of $A$. So, as in the proof of proposition 4, we see that $cl C\subseteq B$. This completes the proof.
\begin{snd}
Let $A$ be a locally convex locally compact subset subset of a locally convex space. Suppose further that $A$ is complete. Let $A_0$ be a compact convex subset of $A$, and $B$ be an open convex subset of $A$ such that $A_0\subseteq B$. Then, there is a continuous function $f:A\to [0,1]$ with compact support such that $f(A_0)=\{1\},\ f(A\backslash B)=\{0\}$, and $1-f$ is a quasi convex function.
\end{snd}
{\textsf Proof:} Find convex open sets $U_0$ and $U_1$ such that $cl(U_0)$ and $cl(U_1)$ are compact and $A\subseteq U_0\subseteq cl(U_0)\subseteq U_1\subseteq cl(U_1)\subseteq B$. As in the proof of first Urysohn lemma, we can find a function  $g:A\to [0,1]$ such that $g(A_0)=\{0\},\ g(A-U_1)=\{1\},\ g$ is a quasi convex and $g$ is continuous. The required function $f$ is $1-g$.

\end{document}